\documentclass[12]{svjour3}

\smartqed  

\usepackage{amssymb, amsmath, enumerate}
\usepackage{graphicx, subfigure}
\usepackage{algorithm, algorithmic}

\newcommand{\bx}{\boldsymbol{x}}
\newcommand{\re}{{\mathrm{re}}}
\newcommand{\im}{{\mathrm{im}}}
\newcommand{\idle}{{\mathrm{idle}}}
\newcommand{\mR}{\mathbb{R}}
\newcommand{\Sov}{\mathbf{H}^1_0(\Omega)}

\DeclareMathOperator{\inte}{int}
\DeclareMathOperator{\exte}{ext}

\newcommand{\closeProof}{ \hfill $\square$  \end{proof} }

\begin{document}

\title{Illusory Shapes via Phase Transition}

\author{Yoon Mo Jung  \and Jianhong Jackie Shen}

\institute{
	Yoon Mo Jung \at
		Department of Computational Science and Engineering, Yonsei University, Korea \\
			\email{ ymjung@yonsei.ac.kr }
		\and
	Jackie Shen (corresponding author) \at
		Department of Industrial and Systems Engineering, University of Illinois, IL 61801, USA \\
		Tel.:  (646) 725-0065 \\
		Fax:  (217) 244-5705 \\
	\email{ jhshen@illinois.edu}
}

\date{\em \small Dedicated to Gil Strang for His 80th Birthday}


\maketitle

\begin{abstract}
We propose a new variational illusory shape (VIS) model via phase fields and phase transitions. It is inspired by the first-order variational illusory contour (VIC) model proposed by Jung and Shen [{\em J. Visual Comm. Image Repres.}, {\bf 19}:42-55, 2008].
Under the new VIS model, illusory shapes are represented by phase values close to 1 while the rest by values close to 0. The 0-1 transition is achieved by an elliptic energy with a double-well potential, as in the theory of $\Gamma$-convergence. The VIS model is non-convex, with the zero field as its trivial global optimum. To seek visually meaningful local optima that can induce  illusory shapes, an iterative algorithm is designed and its convergence behavior is closely studied. Several generic numerical examples confirm the versatility of the model and the algorithm.
\keywords{ Illusory Shapes \and Phase Transition \and Null Hypothesis \and Convergence}
\end{abstract}

\section{Introduction}\label{sec:intro}

The intriguing phenomenon of {\em illusory contours and shapes} (ICS) has been actively pursued in contemporary brain science, neurophysiology, neuroimaging, neural networks, computational neuroscience, and psychophysics~\cite{vision:fMRI2IC-LOC-stanley03,vision:IC-neurons-Lee01,vision:nature-v1forIC-grosof93,vision:neuralnet-interp-fukushima10,vision:nonlinearfeedback-leveille10,vision:salient2IC-yoshino06,vision_Lee_Mumford03}. We first quote  from the recent review paper by Murray and Herrmann~\cite{vision:review-IC-murray13}(2013) to highlight its importance:
\begin{quote} \em
  ``The progressive development of these models ({\em for ICS}) and the data supporting them in many regards highlight what may be considered a general tendency in neuroscience {\em over the past 30 or so years}, namely, progression from ... individual neurons ... to the discrete localization of brain functions ...''
\end{quote}
That is, the persistent scientific interest in ICS has been representative in modern neuroscience, and has profoundly enhanced the {\em systemic} approach to studying complex vision tasks. ICS is fundamentally a system-level vision phenomenon, as contrast to simpler
tasks (e.g., detection of edgelets) by individual simple or complex neurons.

Kanizsa's Triangle illustrated in Fig.~\ref{Fig:Kaniza} is a classic example for ICS. The white triangle in the middle pops out despite that there exist no {\em modal} borders against the white background. Neuronal measurements in the visual cortices of cats or monkeys revealed direct firing responses along such modal edges~\cite{illu_heydt84,illu_heydt89}. This is the {\em very}  nature of illusion, namely, the remarkable ability in reconstruction and re-organization from given visual information with internal structures.

In terms of systems and control,  such ICS phenomena have provided the ideal class of input signals that help identify the complex structures and functionalities of the primate vision system, including, e.g., the bottom-up and top-down schema~\cite{vision_Lee_Mumford03,vision_Wu_Zhu11,vision_Han_Zhu09}. Treating the multi-tier visual cortices and pathways as systems is fundamental. In particular, it implies that models have to be based on {\em systemic} estimations and decisions, instead of isolated neuronal firing patterns. Mathematically, this naturally invites Bayesian decision making over networks, or system-wide optimization~\cite{book_vision_bayes_knill}.

Deterministic variational optimization could be considered as the low temperature limit of the Bayesian decision framework. As in statistical mechanics~\cite{book_mechanics_chandler}, under the low-temperature limit, geometric regularities emerge from statistical gaseous patterns (corresponding to the textural patterns in imaging and vision)~\cite{bImage_ChanShen}. The variational PDE methodology is particularly powerful in handling geometric regularities such as the distance or the curvature~\cite{bImage_AubertKorn,bImage_ChanShen,aubves,luminita_BV}. For ICS, there have been several  variational PDE models in the literature~\cite{illu_sarmalset02,illu_chazhu_levelset03,illu_zhuchan_shape05,illu_jungshen}, as partially reviewed in our earlier work~\cite{illu_jungshen}.

In~\cite{illu_jungshen}, based on the decomposition into {\em real} and {\em imaginary} components, we have proposed a first-order variational illusory contour (VIC) model, which is implemented by the supervised level-set method of Osher and Sethian~\cite{oshset}. Among all the variational-PDE models for ICS that the authors are aware of, this VIC model has the lowest complexity and thus allows detailed analysis of the illusory shapes (as local optima). This will be briefly recapped in Section~\ref{sec:contourmodel}.

The current work has been closely inspired by~\cite{illu_jungshen} and seeks a variational illusory shape (VIS) model based on phase transitions. An illusory shape is ideally represented by the phase value 1.0 and the rest by 0.0. The machinery of $\Gamma$-convergence and phase transitions~\cite{gamma_AmbrosioTor90,gamma_AmbrosioTor92} allows us to construct a phase-field energy that is closely related to the VIC model introduced in~\cite{illu_jungshen}. The new VIS model is more self-contained from modeling to computing. More details are explained in Sections~\ref{sec:represent} and~\ref{sec:new_model}.

Like most phase-field models~\cite{gamma_March92,gamma_MarchDozio97,gamma_ShenPCMS05,gamma_ShenSoftMS06}, the proposed nonlinear VIS model is not convex. The zero field is the global optima but uninteresting. Visually meaningful local optima are hence sought after via an iterative algorithm. The design and analysis of the algorithm will be elaborated in Section~\ref{sec:algorithm}.

In Section~\ref{sec:num_examples}, we summarize the entire algorithm via a pseudo code, and present several generic numerical examples to illustrate the versatility of the model and the algorithm. The work is finally concluded in Section~\ref{sec:conclusion}, where the limitations of the model are also addressed.

\section{Phase Field Representation of Illusory Shapes}
\label{sec:represent}
\begin{figure}[!ht]
 \centering{
 \includegraphics[width=3cm]{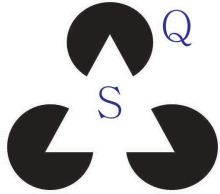}
 \caption{\textit{Kanizsa's Triangle}. The white illusory triangle $S$ in the middle is induced by the black configuration $Q$ consisting of three pac-man disks.}
 \label{Fig:Kaniza} }
\end{figure}

Given a configuration $Q$ in a visual field $\Omega$,
there may emerge an illusory shape $S$ outside $Q$: $S\subseteq \Omega \setminus Q$.
Fig.~\ref{Fig:Kaniza} shows the classic example of \textit{Kanizsa's Triangle}, of which $Q$ consists of three pac-man disks and $S$ refers to the white triangle encompassed and  induced by them.
Mathematically, we assume that $\Omega$ is a bounded Lipschitz convex open domain in $\mathbb{R}^2$, e.g., a Cartesian rectangle or a disk, and that $Q \subset \Omega$ is a {\em compact} Lipschitz sub-domain.

Define the ideal binary phase field $z(\bx)$, with $\bx =(x_1, x_2) \in \Omega$:
\begin{equation*}
z(\bx) =
\left\{
\begin{array}{cl}
  1, & \bx \in S, \\
  0, &  \bx \in \Omega \setminus S.
\end{array}
\right.
\end{equation*}
If $S$ is given, $z(\bx)$ is exactly its indicator function,
and a pixel $\bx \in \Omega$ belongs to $S$ if and only if its phase value is $1$.

In reality, neither $S$ nor $z(\bx)$ is known \textit{a priori}. The phase field approach
attempts to  reconstruct a nontrivial phase field $z(\bx)$'s that can induce illusory shapes. Since binary fields are hard to work with
both theoretically and computationally, one resorts to {\em continuous} fields that are the mollified versions of the binary ones.  They complete the $0$-$1$ phase transition along a tubular neighborhood of the desirable sharp boundaries. Such continuous fields are thus called the {\em diffusive} phase fields in the literature.

For the current task, we impose the following two natural conditions:
\begin{enumerate}[(a)]
 \item $z|_Q \equiv 0$, i.e., illusory shapes can only emerge outside $Q$.

 \item $z|_{\partial \Omega} \equiv 0$. This is because, for variational illusion models like~\cite{illu_jungshen}, one can prove that illusory  shapes can only emerge {\em within} the convex hull of $Q$, and thus also within the interior of $\Omega$.
\end{enumerate}
The main challenge then becomes how to encourage the phase $1.0$ (or close to 1.0) to emerge for configurations that can induce illusory shapes.

\section{Variational Contour Model of Jung and Shen}
\label{sec:contourmodel}

In this section, we briefly review the first-order variational illusory contour (VIC) model by Jung and Shen~\cite{illu_jungshen}. In the next section, the new phase-field model is to be built upon this VIC model.

For simplicity, below we assume all curves to be  piecewise smooth.
Let $\Gamma$ be any simple closed (Jordan) curve on $\Omega \setminus Q^\circ$.
In~\cite{illu_jungshen}, it is decomposed to the {\em real} and {\em imaginary} parts by:
\begin{equation*}
\Gamma_{\re} = \Gamma \cap \partial Q, \quad \Gamma_{\im} = \Gamma \setminus \partial Q.
\end{equation*}
If $\Gamma$ is a genuine illusory contour such as the three sides of Kanizsa's Triangle,
the imaginary part $\Gamma_{\im}$ corresponds to the illusory interpolant.
In Figure \ref{Fig:Kaniza}, for instance, $\Gamma_{\re}$ corresponds to the three corner turns along the disk inducers, and $\Gamma_{\im}$ to the three illusory {\em modal} segments.

For any such a curve $\Gamma$, Jung and Shen proposed the following variational contour model based on the decomposed ``energy''~\cite{illu_jungshen}:
\begin{equation}\label{Eold1}
 E[\Gamma] = E[\Gamma|a,\,b] = a\int_{\Gamma_{\re}}ds + b\int_{\Gamma_{\im}}ds,
\end{equation}
for any pair of properly defined weights $a$ and $b$ with $0<a<b$.
Illusory contours are then defined as the \emph{local} minima.

Compared with several other models~\cite{illu_sarmalset02,illu_chazhu_levelset03,illu_zhuchan_shape05}, the {VIC} model~\eqref{Eold1} is simple and yet powerful enough, as demonstrated by the generic  examples in~\cite{illu_jungshen}. Analytically, the authors were able to establish detailed geometric properties, among which the most important is Theorem 2.11 in~\cite{illu_jungshen}.

\begin{theorem}[Characterization of a Local Minimum] \label{T1:characterization}
 Let $Q$ be a generic compact configuration on an open Lipschitz domain $\Omega$. Assume that $\partial Q$ is piece-wise smooth with finitely many corners (but no kinks).  Let  $\theta_{\min{}},\,
 \theta^\ast_{\min{}} >0$ denote its minimum outer and inner spans. Suppose a given simple closed Jordan curve
 $\Gamma \in \mathcal{C}$ satisfies the following structural conditions:
 \begin{enumerate}[(i)]
 \item (Imaginary Behavior) each connected component $\gamma_\im$ of $\Gamma_\im$ is a
        straight line segment, and no two distinct components share a common hinge;
 \item (Junction Behavior) at any junction point $z \in J[\Gamma]$, the turn $\phi_z <
         \pi/2$, and the idle angle $\phi^\idle_z \ge \pi/2$.
 \end{enumerate}
 Let $\phi_{\max{}}$ denote the maximum turn on $J[\Gamma]$. Then there exists a critical ratio
 $r_c=r_c(\theta_{\min{}}, \phi_{\max{}})< 1$, such that for any $\alpha$ and $\beta$ with
 $r=\alpha/\beta < r_c$, $\Gamma$ is a local minimum to the energy $E[\cdot \mid \alpha,\, \beta]$.
\end{theorem}

We refer to~\cite{illu_jungshen} for the definitions of spans, hinges, turns, and idle angles. Such detailed characterizations are much harder to establish for more complex illusion models~\cite{illu_sarmalset02,illu_chazhu_levelset03,illu_zhuchan_shape05} .

\begin{figure}[!ht]
 \centering{
 \includegraphics[width=7cm]{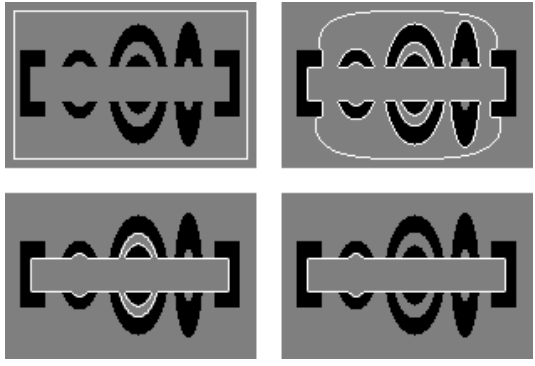}
 \caption{An illusory contour emerges from a complex configuration $Q$ of objects, as captured by the supervised level-set method of Osher and Sethian~\cite{oshset} for the VIC model in  Eqn.~\eqref{Eold1} proposed by Jung and Shen~\cite{illu_jungshen}. }
 \label{Fig1}
 }
\end{figure}

In~\cite{illu_jungshen},  illusory shapes (as the local minima to $E[\Gamma|a,\,b]$)
have been computed via a supervised level-set scheme of Osher and Sethian~\cite{oshset},
as applied to the approximate energy (Eqn.~(7) in~\cite{illu_jungshen}):
\begin{equation}
E_\sigma[\Gamma|\alpha,\, \beta] = \alpha\int_{\Gamma}ds + \beta\int_{\Gamma} g(\lvert \nabla \chi_{Q, \sigma} (s)\rvert)\, ds,
\end{equation}
with $\alpha = a + b$, $\beta = b$ for the targeted energy $E[\Gamma|a,\, b]$ in Eqn.~\eqref{Eold1}. Here $\sigma \ll 1$ denotes a small diffusion or mollification scale, so that the binary indicator $\chi_Q$ is mollified to a smooth approximation $\chi_{Q, \sigma}$.
The function $g$ could be any positive function satisfying: $g(0^+) = 1$, $g(+\infty) = 0$, for example, $g(p) =\exp{(-p^2)}$ or
$g(p) = 1/(1+p^2)$. It was shown in~\cite{illu_jungshen} that for any admissible contour $\Gamma$,
$E_\sigma[\Gamma|\alpha,\, \beta] \rightarrow E[\Gamma|a,\, b]$ as $\sigma \to 0^+$. The critical ratio $r_c$ for $a/b$ in Theorem \ref{T1:characterization} now transfers to  $R_c = r_c/(1-r_c)$ for $\alpha/\beta$, accordingly.

\section{A Phase Transition Model for Illusory Shapes}
\label{sec:new_model}
In this section, we develop a new phase transition model for illusory shapes, based on the discussion in the preceding two sections.

We first define the inducing scalar field via $G(\bx) = \alpha +\beta g(\lvert \nabla \chi_{Q, \sigma} (\bx)\rvert)$, for any $\bx\in \Omega$.
Alternatively, one could employ any optimal phase field $z^{\text{MS}}$ from the phase transition approximation to the Mumford-Shah model~\cite{gamma_AmbrosioTor90,gamma_AmbrosioTor92}:
\begin{equation*}
G(\bx) = \alpha + \beta z^{\text{MS}}_\sigma(\bx), \quad \bx \in \Omega.
\end{equation*}
The basic requirements for $z^{\text{MS}}_\sigma(\bx)$ are:
(i) the phase transition bandwidth is small: $\sigma \ll 1$, (ii) $u_0(\bx) = \chi_{Q}(\bx)$  is used as the input to the Mumford-Shah model, and (iii) $z^{\text{MS}}_\sigma(\bx)\in [0, \, 1]$ with $z^{\text{MS}}_\sigma \simeq 0$ along $\partial Q$ and $\simeq 1$  away from $\partial Q$.

In either way, one has
\begin{equation*}
G(\bx) \simeq
\left\{
\begin{array}{ll}
\alpha, &  \bx \,\, \text{near or along}\,\, \partial Q,\\
\alpha + \beta, & \text{otherwise}.
\end{array}
\right.
\end{equation*}
We call $G$ a   \emph{canyon} function associated with the given $Q$, since near $\partial Q$  the values of $G$ quickly drop from $\alpha + \beta$ to $\alpha$. For the canyon effect to be more salient, it is natural to require the dropping size $\beta \gg \alpha$, which is also in accordance with the critical ratio requirement $\alpha/\beta <R_c$ in the preceding section.

Let $\Sov$ denote the Sobolev space of functions on $\Omega$ with zero boundary traces.
It hosts all diffusive phase fields discussed in Section~\ref{sec:represent}.
For any phase field $z\in \Sov$, we introduce the following energy:
\begin{equation}\label{model:phase}
E_\epsilon[z] = \int_\Omega \left[\frac{\epsilon}{2} \lvert\nabla z \rvert^2+ \frac{(1-z)^2z^2}{2\epsilon} \right]\cdot G\,d\bx
 + \lambda \int_\Omega \chi_Q \frac{z^2}{2\epsilon}\,d\bx,
\end{equation}
where the positive weight $\lambda$ is in the same order as  of $G$.
When $\alpha < \beta = O(1)$ so that $G=O(1)$, we simply set $\lambda =1$ in all our computational examples later on. The small ``diffusivity'' parameter $\epsilon$ defines the intended transition bandwidth between 0.0 and 1.0.

The two terms have been motivated  as follows. The second term is:
\begin{equation*}
\lambda \int_\Omega \chi_Q \frac{z^2}{2\epsilon}\,d\bx = \lambda \int_Q \frac{z^2}{2\epsilon}\,d\bx.
\end{equation*}
Thus for $\epsilon \ll 1$, it acts as a soft way to enforce $z  = 0$ on $Q$, which is the condition (a) imposed in Section~\ref{sec:represent}.  {Also notice that the condition (b) is met automatically since $z \in \Sov$.}

For the first term in $E_\epsilon[z]$, define the Borel measure $\mu_\epsilon$ for any given $z\in  \Sov$ by:
\begin{equation*}
d\mu_\epsilon = \left[\frac{\epsilon}{2} \lvert\nabla z \rvert^2+ \frac{(1-z)^2z^2}{2\epsilon}\right] \, d\bx.
\end{equation*}
For any smooth Jordan curve $\Gamma$ with arclength element $ds$,
we show semi-heuristically that $d\mu_\epsilon \rightarrow \frac{1}{6}ds$ as $\epsilon \rightarrow 0^+$, or more generally via the 1-dimensional Hausdorff measure $d\mathcal{H}^1_\Gamma$ that $d\mu_\epsilon \rightarrow \frac{1}{6}d\mathcal{H}^1_\Gamma$  (with $\mathcal{H}^1_\Gamma(A) = \mathcal{H}^1(\Gamma \cap A)$) for a properly
designed sequence of phase fields {$\{z_\epsilon\}_\epsilon$}.
Assume that $\Gamma$ is the $0$-level set of a smooth function $f(\bx)$: $\Gamma = f^{-1}(0)$, and that $f$ is \emph{regular} along $\Gamma$ in the sense that $\nabla f(\omega) \neq 0$ for all $\omega \in \Gamma$. With at most a sign flip, one could further assume that  $\nabla f(\omega)$ points towards the inside of $\Gamma$. Then by the Tubular Neighborhood Theorem~\cite{booth86}, there exists some $\delta_0 >0$ such that the map:
$$\varphi: (\omega, n)\in \Gamma \times (-2\delta_0, 2\delta_0)
\quad \rightarrow \quad \bx = \omega + n \cdot \frac{\nabla f}{\lvert\nabla f \rvert} \in \Omega, $$
is a diffeomorphism between $\Gamma \times (-2\delta_0, 2\delta_0)$ and an open tubular neighborhood $B_{2\delta_0}$ of $\Gamma$ in $\Omega$.
Let $S$ denote the logistic function:
$$S(t) = \frac{1}{1+e^{-t}}, \quad t\in \mR^1.$$
For any $\epsilon \ll \delta_0$, we then construct a special phase field $z_\epsilon \in \Sov$:
\[
z_\epsilon(\bx) =
\left\{
\begin{array}{ll}
S\left(\frac{n}{\epsilon}\right), & \bx = \varphi(\omega, n)\in B_{\delta_0},\\
S\left(\frac{\delta_0}{\epsilon}\right), &  \bx \notin  B_{\delta_0} \,\, \text{and} \,\, \bx \in \inte(\Gamma),\\
S\left(-\frac{\delta_0}{\epsilon}\right), &   \bx\notin B_{\delta_0} \,\, \text{and} \,\, \bx \in \exte(\Gamma).
\end{array}
\right.
\]
Here $\inte(\Gamma)$ and $\exte(\Gamma)$ denote the interior and exterior domains of $\Gamma$, which are well defined according to the renowned Jordan Curve Theorem in topology~\cite{topo_thmJordan_hales2007}.

One can show that for any \emph{continuous} function $\phi$ on $\Omega$:
\begin{equation}
\int_\Omega \phi\,d\mu_\epsilon \rightarrow \frac{1}{6}\int_\Omega \phi \,d\mathcal{H}^1_\Gamma
= \frac{1}{6}\int_\Gamma \phi \,ds, \quad \text{as}\quad \epsilon \rightarrow 0^+,
\end{equation}
where the multiplier $\frac{1}{6}$ is the total variation of $\displaystyle F(z) = \frac{z^2}{2} -  \frac{z^3}{3}$ from $z = 0$ to $z=1$, and $F(z)$ is the primitive of $z(1-z)$ so that $F'(z)=z(1-z)$. As a result, for example, if $\Gamma$ represents the illusory Kanizsa Triangle in Fig.~\ref{Fig:Kaniza} and $\phi=G$, as $\varepsilon \to 0$,  one has
\[
6\int_\Omega G\,d\mu_\epsilon \simeq a\int_{\Gamma_{\re}}ds + b\int_{\Gamma_{\im}}ds,
\quad a=\alpha, \; b= \alpha + \beta.
\]
This gives a phase-field interpretation of the model by Jung and Shen~\cite{illu_jungshen}.
A similar argument could also be found in~\cite{gamma_ShenPCMS05}, for example.
But rigorous treatment is only offered by the theory of $\Gamma$-convergence
approximation~\cite{bGamma_Braides02,bGamma_DalMaso92}.

In summary, the phase transition model $E_\epsilon[z]$ for $z\in \Sov$ proposed in Eqn.~\eqref{model:phase} generalizes the contour model $E[\Gamma]$ of Jung and Shen~\cite{illu_jungshen} in Eqn.~\eqref{Eold1} to continuous phase fields. It also implements the two conditions (a) and (b) imposed in Section~\ref{sec:represent}.

One notices that $z_* \equiv 0$ is the global minimum of $E_\epsilon[z]$ but uninteresting. We thus define any local minimum $z_\epsilon \neq z_*$ to be an \emph{illusory phase field}, and
\begin{equation}\label{IlluShape}
S_\epsilon = \{ \bx\in \Omega: z_\epsilon > 1/2\}
\end{equation}
to be the associated \emph{illusory shape}.

\begin{theorem}\label{Thm2}
Suppose $z_\epsilon$ is a local minimum of $E_\epsilon[z]$ and is not always $0$ on $\Omega$.
Then the associated illusory shape $S_\epsilon \neq \emptyset$.
\end{theorem}
\begin{proof}
Otherwise assume $S_\epsilon = \emptyset$. Then $z_\epsilon \le 1/2$ for any $\bx \in \Omega$.
For any $t\in(0, \,1)$, define $z^{(t)}(\bx) = tz_\epsilon$. Then $z^{(t)} \in \Sov$,
and for $t$ close to $1$,  $z^{(t)}$ is a perturbation to $z^{(1)} = z_\epsilon$. One has:
\begin{equation}\label{eqn:perturb}
E_\epsilon[z^{(t)}] = t^2\left[\int_\Omega \frac{\epsilon}{2} \lvert\nabla z_\epsilon \rvert^2 \cdot G\,d\bx
 +\lambda \int_\Omega \chi_Q \frac{z_\epsilon^2}{2\epsilon}\,d\bx   \right]
 +\int_\Omega \frac{\Phi(tz_\epsilon)}{2\epsilon} \cdot   G\,d\bx,
\end{equation}
where $\Phi(z) = (1-z)^2z^2$ is the double-well potential.

Since $\Phi(z)$ is decreasing on $[-\infty,\, 0]$ and increasing on $(0, \, 1/2]$, and $z_\epsilon(\bx) \le 1/2$
for any $\bx \in \Omega$, one has:
$$  \Phi( tz_\epsilon(\bx)) \le  \Phi(z_\epsilon(\bx)), \quad \text{for any } \bx \in \Omega \text{ and } t\in(0,\, 1).$$
Thus \eqref{eqn:perturb} implies that for any $t\in(0,\, 1)$, $E_\epsilon[z^{(t)}] \le E_\epsilon[z_\epsilon]$. In fact, one must have $E_\epsilon[z^{(t)}] < E_\epsilon[z_\epsilon]$. This is because
$$ \int_\Omega \frac{\epsilon}{2} |\nabla z_\epsilon |^2 \cdot G\,d\bx
\ge \frac{\epsilon\alpha}{2} \int_\Omega  \lvert\nabla z_\epsilon \rvert^2 \,d\bx
= \frac{\epsilon\alpha}{2} |z_\epsilon|^2_1,$$
with $\alpha =\min G > 0$. Notice that in $\Sov$, the first-order homogeneous semi-norm $|\cdot|_1$ is actually a norm due to the zero trace~\cite{adafou}. Since it has been assumed that $z_\epsilon \neq 0 $ in $\Sov$, one must have $|z_\epsilon|^2_1> 0$.
Thus the first term alone in \eqref{eqn:perturb} shows $E_\epsilon[z^{(t)}] < E_\epsilon[z_\epsilon]$ for $t\in(0,\, 1)$. This contradicts to the assumption that $z_\epsilon$ is a local minimum since $z^{(t)}$'s are its small perturbations for $t \simeq 1$.
\closeProof

\section{Null Hypotheses and An Iterative Algorithm}
\label{sec:algorithm}

We now design a specific iterative algorithm to seek visually meaningful local minima of the proposed energy:
\begin{equation}\label{model:phase2}
E_\epsilon[z] = \int_\Omega \left[\frac{\epsilon}{2} \lvert\nabla z \rvert^2+ \frac{(1-z)^2z^2}{2\epsilon} \right]\cdot G\,d\bx
 + \lambda \int_\Omega \chi_Q \frac{z^2}{2\epsilon}\,d\bx,
\end{equation}
for $z\in\Sov$, given $G \in \mathcal{C}(\Omega)$ and $G \ge 0$.
For the current work, $G$ is specifically constructed as in the opening of the preceding section. In particular, we assume that
\begin{equation}\label{condG}
\alpha = \min G > 0\quad \text{and} \quad \alpha + \beta = \max G,
 \text{ with } \beta \ge \alpha.
\end{equation}
For optimization, $\lambda$ could be ``absorbed'' into $G$ via $G/\lambda$.
Thus we assume $\beta = O(1)$, and $\lambda = 1$.

In the literature of non-convex optimization, there has been much discussion on finding the global optima (e.g., stochastic or deterministic annealing~\cite{comput_annealing_kirkpatrick1983}). But no universal methodologies exist for locating the local optima that are of practical interest. Generally it has to be problem specific, and the two key components are:
\begin{enumerate}[(i)]
\item the choice of an initial guess or state, and
\item the design of an iterative searching algorithm.
\end{enumerate}

For the initial guess, as motivated by hypothesis testing in statistics~\cite{book_stat_elem_freedman}, we work with the following \emph{null hypothesis}:
\begin{equation}\label{hypothesis}
\text{``There indeed exists an illusory shape somewhere outside $Q$."}
\end{equation}
As the information of the illusory shape is unknown {\em a priori},
we start with following initial guess:
\begin{equation}\label{eqn:iniguess}
z_0(\bx) = 1\cdot(1-\chi_Q) + 0\cdot \chi_Q = 1-\chi_Q.
\end{equation}
It conservatively assigns phase $1$ to all pixels outside the given configuration $Q$. Strictly speaking, $z_0\notin \Sov$. This is not an issue since $z_0$ is to be used in an iterative algorithm: $z_n \rightarrow z_{n+1}$, and $z_n$'s generated afterwords all belong to $\Sov$ for $n = 1, 2, \cdots$.
Alternatively, one may apply diffusion to $z_0:\,\, u_t =\Delta u, \,\, u|_{\partial \Omega} = 0, \,\, u|_{t=0} = z_0$,
and use $u(\cdot, \delta)\in \Sov$ instead as the initial guess for some $\delta \ll 1$.

Next, to design the iterative algorithm, we first compute the Euler-Lagrange equation of $E_\epsilon$ in Eqn.~\eqref{model:phase2}:
$$-\nabla\cdot(\epsilon G \nabla)z +\frac{G}{\epsilon}(z-3z^2+2z^3) + \frac{\lambda \chi_Q}{\epsilon}\cdot z = 0,$$
which can be rearranged as:
\begin{equation}\label{eqn:pde}
-\nabla\cdot(\epsilon^2 G \nabla)z + (G(1+2z^2)+\lambda \chi_Q)z = 3Gz^2.
\end{equation}
This is a nonlinear elliptic equation on $\Omega$ with boundary condition $z|_{\partial \Omega} = 0$.
As discussed earlier, the global minimum  $z_\ast\equiv 0$ is a solution.
We are interested in the non-zero solutions that are associated with the local minima of $E_\epsilon[z]$.

Our proposed iterative algorithm is to solve the following linear elliptic equation for $z_{n+1}$, given a current guess $z_n$:
\begin{equation}\label{eqn:linpde}
\begin{aligned}
&-\nabla\cdot(\epsilon^2 G \nabla)z + g_nz = f_n, \quad z|_{\partial \Omega} = 0, \quad \text{with}\\
& g_n = G(1+2z^2_n)+\lambda \chi_Q,\quad \mbox{and}\; f_n=3Gz^2_n.
\end{aligned}
\end{equation}

\begin{theorem}
Suppose that $G$ satisfies Eqn.~\eqref{condG} and $z_n \in \mathbf{L}^\infty(\Omega)$.
Then a weak solution $z = z_{n+1}\in \Sov$ to Eqn.~\eqref{eqn:linpde} always exists and is unique in the sense that, for any $u \in \Sov$,
\begin{equation}\label{weakform}
(\epsilon^2 G\nabla z_{n+1}, \, \nabla u) + (g_nz_{n+1}, \, u) = (f_n,\, u),
\end{equation}
where $(\cdot,\, \cdot)$ denotes the canonical inner product in $\mathbf{L}^2(\Omega,\, \mR^k)$.
\end{theorem}
\begin{proof}
In $\Sov$, define the symmetric bilinear function $\langle \cdot,\, \cdot \rangle_1$ via:
$$\langle u,\, v \rangle_1  = (\epsilon^2 G\nabla u, \, \nabla v) + (g_n u, \, v),\quad u,v\in \Sov,$$
and denote the canonical inner product in $\Sov$ by $\langle \cdot,\, \cdot \rangle_0 $, which is defined by:
$$\langle u,\, v \rangle_0  = (\nabla u, \, \nabla v) + (u, \, v).$$
By Eqn.~\eqref{condG},
$$\epsilon^2\alpha (\nabla u, \, \nabla u) \le (\epsilon^2 G\nabla u, \, \nabla u) \le \epsilon^2(\alpha+\beta)(\nabla u, \, \nabla u),$$
and by the definition of $g_n$ in~\eqref{eqn:linpde},
$$\alpha (u, \, u) \le  (G u, \, u)  \le  (g_nu, \, u) \le \big((\alpha+\beta)(1+2\|z_n\|^2_\infty)+\lambda\big)  (u, \, u). $$
Therefore, $\langle \cdot,\, \cdot \rangle_1$ is an inner product equivalent to $\langle u,\, v \rangle_0$. Noticing that
$$|(f_n,\, u)| \le 3(\alpha+\beta)\|z_n\|^2_\infty(1,\, |u|) \le 3(\alpha+\beta)\|z_n\|^2_\infty\sqrt{|\Omega|}\sqrt{\langle u, \,u \rangle_0}, $$
$(f_n,\, \cdot)$ must be a continuous linear function for $\langle \cdot ,\, \cdot \rangle_0$ and thus also $\langle \cdot ,\, \cdot \rangle_1$.
Applying Riesz representation theorem~\cite{fol} to $\langle \cdot , \cdot \rangle_1$ and $(f_n,\, \cdot)$,
one concludes that there exists a unique $z_{n+1} \in \Sov$, such that
$$\langle z_{n+1}, \, u \rangle_1 \equiv (f,\, u),\quad \text{for any } u\in \Sov,$$
which is precisely Eqn.~\eqref{weakform}.
\closeProof

We also have an energy-form description for $z_{n+1}$ under given $z_n$.

\begin{proposition}\label{Prop1}
$z = z_{n+1}\in\Sov$ is the unique weak solution to Eqn.~\eqref{eqn:linpde}
if and only if $z_{n+1} = \arg\min_z E_\epsilon[z|z_n]$, where
\begin{equation}\label{model:phase3}
E_\epsilon[z|z_n] = \int_\Omega \left[\frac{\epsilon}{2} \lvert\nabla z \rvert^2+ (1+2z_n^2)\frac{(z-\gamma_n)^2}{2\epsilon} \right]\cdot G\,d\bx
 + \lambda \int_\Omega \chi_Q \frac{z^2}{2\epsilon}\,d\bx,
\end{equation}
with $\gamma_n = 3z^2_n/(1+2z^2_n)$. Notice that $E_\epsilon[\cdot |z_n]$ is strictly convex.
\end{proposition}
\begin{proof} By definition, for any $z, u \in \Sov$,
\begin{equation} \label{E:expansion}
E_\epsilon[z+u|z_n] = E_\epsilon[z|z_n] + E^*_\epsilon[u|z_n] + J_\epsilon[z, u|z_n],
\end{equation}
with
$$E^*_\epsilon[u|z_n] = \int_\Omega \left[\frac{\epsilon}{2} \lvert\nabla u \rvert^2+ (1+2z_n^2)\frac{u^2}{2\epsilon} \right]\cdot G\,d\bx
 + \lambda \int_\Omega \chi_Q \frac{u^2}{2\epsilon}\,d\bx, $$
and
$$J_\epsilon[z, u|z_n] = (\epsilon G \cdot \nabla z, \, \nabla u) +  \frac{1}{\epsilon}(G(1+2z_n^2)(z-\gamma_n), \, u)
+\frac{\lambda}{\epsilon}( \chi_Q \cdot z, \, u).$$
Since $E^*_\epsilon[u|z_n]\ge 0$, $E^*_\epsilon[tu|z_n]=t^2E^*_\epsilon[u|z_n]$, and
$J_\epsilon[z, tu|z_n] = tJ_\epsilon[z, u|z_n]$, we conclude from Eqn.~\eqref{E:expansion} by using $t \ll 1$ that
$$z_{n+1} = \arg\min_z E_\epsilon[z|z_n]\quad \text{iff}\quad J_\epsilon[z_{n+1}, u|z_n] = 0,
\,\,\text{for any } u\in \Sov.$$
On the other hand,
$$\epsilon \cdot J_\epsilon[z, u|z_n] = (\epsilon^2 G\nabla z, \, \nabla u) +  (g_nz,\, u)-(f_n,\, u),$$
which is precisely Eqn.~\eqref{weakform} with $z=z_{n+1}$.  Thus $z_{n+1}=\mathrm{argmin} E_\epsilon[z | z_n]$ if and only if it satisfies Eqn.~\eqref{weakform}. {Moreover,  $E^*_\epsilon[u|z_n] > 0$ for $u \neq 0$ guarantees strict convexity.} This completes the proof.
\closeProof

This leads to the following desirable property for $z_n$'s to faithfully approximate $0$-$1$ binary phases.
\begin{proposition}\label{Prop2}
Let $z_{n+1}\in \Sov$ be the unique weak solution to the iterative algorithm in Eqn.~\eqref{eqn:linpde}
given $z_n$. If $z_n(\bx) \in [0, 1]$ for any $\bx\in\Omega$, so must be $z_{n+1}$.
\end{proposition}
\begin{proof} Given $z_{n+1} \in \Sov$, define truncation $z^{[0, 1]}_{n+1}$ by
$$
z^{[0, 1]}_{n+1}(\bx) =
\left\{
\begin{array}{ll}
z_{n+1}(\bx), & \text{if} \quad 0 < z_{n+1}< 1,\\
1, & \text{if}\quad z_{n+1}\ge 1,\\
0, &  \text{if}\quad z_{n+1} \le 0.
\end{array}
\right.
$$
Then $z^{[0, 1]}_{n+1}\in \Sov$, and it is well known~\cite{adafou} that truncation does not increase
the norm of a gradient in Sobolev spaces. Then one must have:
$$\int_\Omega \frac{\epsilon}{2} \left\lvert\nabla z^{[0, 1]}_{n+1} \right\rvert^2 G\,d\bx
\le \int_\Omega \frac{\epsilon}{2} \lvert\nabla z_{n+1} \rvert^2G\,d\bx. $$
Since $z_n\in[0,1]$, one has $z_n^2 \le 1$, and for the scalar field $\gamma_{n}$ in Eqn.~\eqref{model:phase3},
$$\gamma_n = \frac{3z_n^2}{1+2z_n^2}\le 1.$$
The one must have for any $\bx\in \Omega$,
$$(z^{[0, 1]}_{n+1}(\bx) - \gamma_n )^2 \le (z_{n+1}(\bx) - \gamma_n )^2\quad \text{and}\quad
\left(z^{[0, 1]}_{n+1}(\bx)\right)^2 \le (z_{n+1}(\bx))^2. $$
Examining the expression in Eqn.~\eqref{model:phase3} thus shows:
$$E_\epsilon\left[z^{[0, 1]}_{n+1}\big|z_n\right]\le E_\epsilon[z_{n+1}|z_n]. $$
Due to the uniqueness result from the preceding proposition,
$$z_{n+1}(\bx)\equiv z^{[0, 1]}_{n+1}(\bx)\in [0, 1], \quad \bx\in\Omega.$$
This completes the proof.
\closeProof

The next theorem reveals that the phase field sequence $(z_n)$ defined by the iterative algorithm in~\eqref{eqn:linpde} is indeed an energy decreasing sequence.

\begin{theorem}\label{thm4}
Let $\{z_n\}^\infty_{n=1}$ be the sequence of phase fields generated by the algorithm in Eqn.~\eqref{eqn:linpde}, started from the null hypothesis $z_0$ in Eqn.~\eqref{eqn:iniguess}. Then
$$E_\epsilon[z_1]\ge E_\epsilon[z_2]\ge \cdots \ge E_\epsilon[z_n] \ge \cdots $$
Furthermore, for $n \ge 1$,
\begin{equation}\label{eqn:difference}
E_\epsilon[z_n]- E_\epsilon[z_{n+1}] \ge \int_\Omega \frac{G}{2\epsilon}
(z_{n+1} - z_n)^2\left(2z_{n}+4z_{n+\frac{1}{2}}\left(1-z_{n+\frac{1}{2}}\right)\right)d\bx,
\end{equation}
where $z_{n+\frac{1}{2}} =  (z_n + z_{n+1})/2$.
\end{theorem}
\begin{proof}
By the null hypothesis $z_0$ and Proposition~\ref{Prop2}, we have
$$0 \le z_n(\bx) \le 1, \quad \text{for all } \bx\in\Omega \quad\text{and}\quad n = 1, 2, \cdots.$$
Then $z_{n+\frac{1}{2}}(\bx) =  (z_n(\bx) + z_{n+1}(\bx))/2\in [0, 1]$ for any $\bx \in \Omega$.
Thus Eqn.~\eqref{eqn:difference} implies $E_\epsilon[z_n]\ge E_\epsilon[z_{n+1}]$ for $n\ge 1$ since the
integrand is nonnegative. It suffices to establish~\eqref{eqn:difference} for $n\ge 1$.
By Proposition~\ref{Prop1}, since $z_{n+1} = \arg\min_z E_\epsilon[z|z_n]$, one has
$$ E_\epsilon[z_{n+1}|z_n] \le E_\epsilon[z_n|z_n]. $$
This spells out to be:
$$\begin{aligned}
 &\int_\Omega \frac{\epsilon G}{2} \lvert\nabla z_{n+1} \rvert^2 \,d\bx
 + \lambda \int_\Omega \chi_Q \frac{z_{n+1}^2}{2\epsilon}\,d\bx
 + \int_\Omega \frac{G}{2\epsilon} \left((1+2z_n^2)z_{n+1}^2 - 6z_n^2 z_{n+1}\right) d\bx\\
  \le & \int_\Omega \frac{\epsilon G}{2} \lvert\nabla z_n \rvert^2\, d\bx
 + \lambda \int_\Omega \chi_Q \frac{z_{n}^2}{2\epsilon}\,d\bx
 + \int_\Omega \frac{G}{2\epsilon}  \left((1+2z_n^2)z_{n}^2 - 6z_n^3\right) d\bx.
\end{aligned} $$
Then by the definition of $E_\epsilon[z]$ in Eqn.~\eqref{model:phase2},
$$ \begin{aligned}
& E_\epsilon[z_n]- E_\epsilon[z_{n+1}] \\
 \ge & \int_\Omega \frac{G}{2\epsilon} \left((1+2z_n^2)z_{n+1}^2 - 6z_n^2 z_{n+1}\right) d\bx
- \int_\Omega \frac{G}{2\epsilon} \left((1+2z_n^2)z_{n}^2 - 6z_n^3\right) d\bx \\
&+ \int_\Omega \frac{G}{2\epsilon} z_n^2(1-z_n)^2 \,d\bx
-\int_\Omega \frac{G}{2\epsilon} z_{n+1}^2(1-z_{n+1})^2 \,d\bx.
\end{aligned}$$
Expanding and re-organizing the integrands on the right, one arrives at the inequality~\eqref{eqn:difference}.
\closeProof

\begin{proposition}\label{Prop3}
Following the same assumptions of Theorem \ref{thm4},
we conclude that there must exist some $E_*\ge 0$,
such that $\lim_{n \rightarrow \infty} E_\epsilon[z_n] = E_*$.
\end{proposition}
This is because any bounded monotonic sequence must converge.

\begin{proposition}\label{Prop4}
If  $E_* = \lim_{n \rightarrow \infty} E_\epsilon[z_n] = 0$,  $\{z_n\}^\infty_{n=1}$ must converge to the
global minimum $z_* \equiv 0$ in $\Sov$.
\end{proposition}
\begin{proof}
By definition of $E_\epsilon[\cdot]$ in Eqn.~\eqref{model:phase2}, and the assumptions on $G$ in Eqn.~\eqref{condG},
$$E_\epsilon[z] \ge \frac{\epsilon \alpha}{2} \int_\Omega  |\nabla z|^2 \,d\bx
\ge \frac{\epsilon \alpha}{2}C_\Omega \|z\|^2_{\Sov},$$
where $C_\Omega$ is a positive constant only depending on the domain. The second inequality holds since the trace along $\partial \Omega$ vanishes~\cite{adafou}. Thus $E_\epsilon[z_n]\rightarrow 0$
implies that  $z_n \rightarrow 0$ in $\Sov$.
\closeProof

Let $\rho_n = E_\epsilon[z_n]- E_\epsilon[z_{n+1}]$
denote the successive energy improvement. Then
$$\sum^\infty_{n =1} \rho_n = E_\epsilon[z_1] - E_* < \infty.$$
The series $(\rho_n)$ are said to converge in a generalized quadratic power-law (GQPL) if
$$\sum^\infty_{n =1}\sqrt{\rho_n} < \infty.$$
For instance, $\displaystyle \rho_n = 1/(1+n)^{2+\delta}$ or $1/ ((1+n)^2 \log^{2+\delta}(1+n))$ for any $\delta > 0$.

\begin{proposition}
Suppose $\{E_\epsilon[z_n]\}^\infty_{n=1}$ converges to some $E_* > 0$ in GQPL.
Suppose there exists
a measurable set $S \subseteq\Omega$ with a positive Lebesgue measure,
some constant $C>0$ and integer $N$ such that for any $n>N$,
$$z_n(\bx) \ge C, \quad \text{for }  \bx \in S. $$
Then $\{z_n|_S\}^\infty_{n=1}$ converges in $\mathbf{L}^2(S)$.
\end{proposition}
\begin{proof}
For $n > N$,
$$2z_n + 4z_{n+\frac{1}{2}}(1-z_{n+\frac{1}{2}}) \ge 2 z_n \ge 2C.$$
Then Eqn.~\eqref{eqn:difference} leads to:
$$\rho_n \ge \int_S \frac{G}{2\epsilon} (z_{n+1}-z_n)^2 (2C)\,d\bx
\ge \frac{\alpha C}{\epsilon}\|z_{n+1}-z_n\|^2_{\mathbf{L}^2(S)}.$$
Thus,
$$\|z_{n+1}-z_n\|_{\mathbf{L}^2(S)} \le \sqrt{\frac{\epsilon}{\alpha C}} \cdot \sqrt{\rho_n}.$$
Since $\sum^\infty_{n =1}\sqrt{\rho_n} < \infty$, it implies that {$\{z_n\}_{n=1}^\infty$} is a Cauchy sequence in $\mathbf{L}^2(S)$.
\closeProof

The result can be intuitively interpreted as follows for the half-way threshold $C = 0.5$.
One starts with the null hypothesis (via $z_0$) that the illusory shape is $S_\epsilon = \Omega\setminus Q$.
As the iteration progresses, some pixels are rejected if $z_n(\bx) < 0.5$.
But suppose there exists a positive set $S$, such that all $z_n$'s after some $N$
consistently vote for it in the sense of $z_n(\bx) > 0.5$. Then the voting must be ``directional'' or converging, and $S$ has to be part of the final illusory shape.

\section{Numerical Implementation and Examples}
\label{sec:num_examples}
In this section, we briefly describe the numerical scheme of the proposed model and
present several generic computational examples.

Algorithm~\ref{algo} below offers a pseudocode block describing the major computational steps for the proposed model. For the initial guess $z_0$, we have adopted the null hypothesis~\eqref{hypothesis}. The core iteration formula for updating $z_{n+1}$ from $z_n$ follows  Eqn.~\eqref{eqn:linpde}.
Convergence analysis of the algorithm has been partially given in the preceding section.

\begin{algorithm}
\caption{Illusory Shapes via Phase Transition} \label{algo}
\begin{algorithmic}
\STATE {\bf input}: the configuration $Q$ represented via its indicator $\chi_Q$;
\STATE {\bf initialize:} $z_0(\bx) = 1\cdot(1-\chi_Q) + 0\cdot \chi_Q = 1-\chi_Q$ as the null hypothesis;
\STATE {\bf pre-process:} the canyon funtion $G(\bx) = \alpha +\beta g(\lvert \nabla \chi_{Q, \sigma} (\bx)\rvert)$ or $\alpha +\beta z^{\mathrm{MS}}_\sigma(\bx)$;
\WHILE{$\Vert z_{n+1} - z_{n} \Vert > \delta $}
\STATE Solve for $z_{n+1}$ from:
\begin{equation*}\label{eqn:algo}
\begin{aligned}
& -\nabla\cdot(\epsilon^2 G \nabla)z + g_nz = f_n, \quad\text{with}\\
& g_n = \,G(1+2z^2_n)+\lambda \chi_Q,\, f_n=3Gz^2_n,\text{ and } z|_{\partial \Omega} = 0.
\end{aligned}
\end{equation*}
\ENDWHILE
\end{algorithmic}
\end{algorithm}

Therefore, the core of the algorithm is an elliptic solver, which can be found in the standard literature of computational PDE's~\cite{book_numerical_Sauer}. Some typical parameters generating the examples herein are given as follows:
\[
	\alpha =0.1;\; \beta =1.0;\; \lambda =1.0;\;\; \mbox{and}\;\; \epsilon = (2\sim 4) \cdot h.
\]
Here $h$ denotes the grid/pixel size and is defined in such a way that the longest size of the image domain $\Omega$ is always normalized to the unit length.
The convergence tolerance is set to be $\delta =10^{-6}$.

\begin{figure}[h!]
  \centering
  \subfigure[Kaninza Triangle]{
   \includegraphics[width=3.5cm]{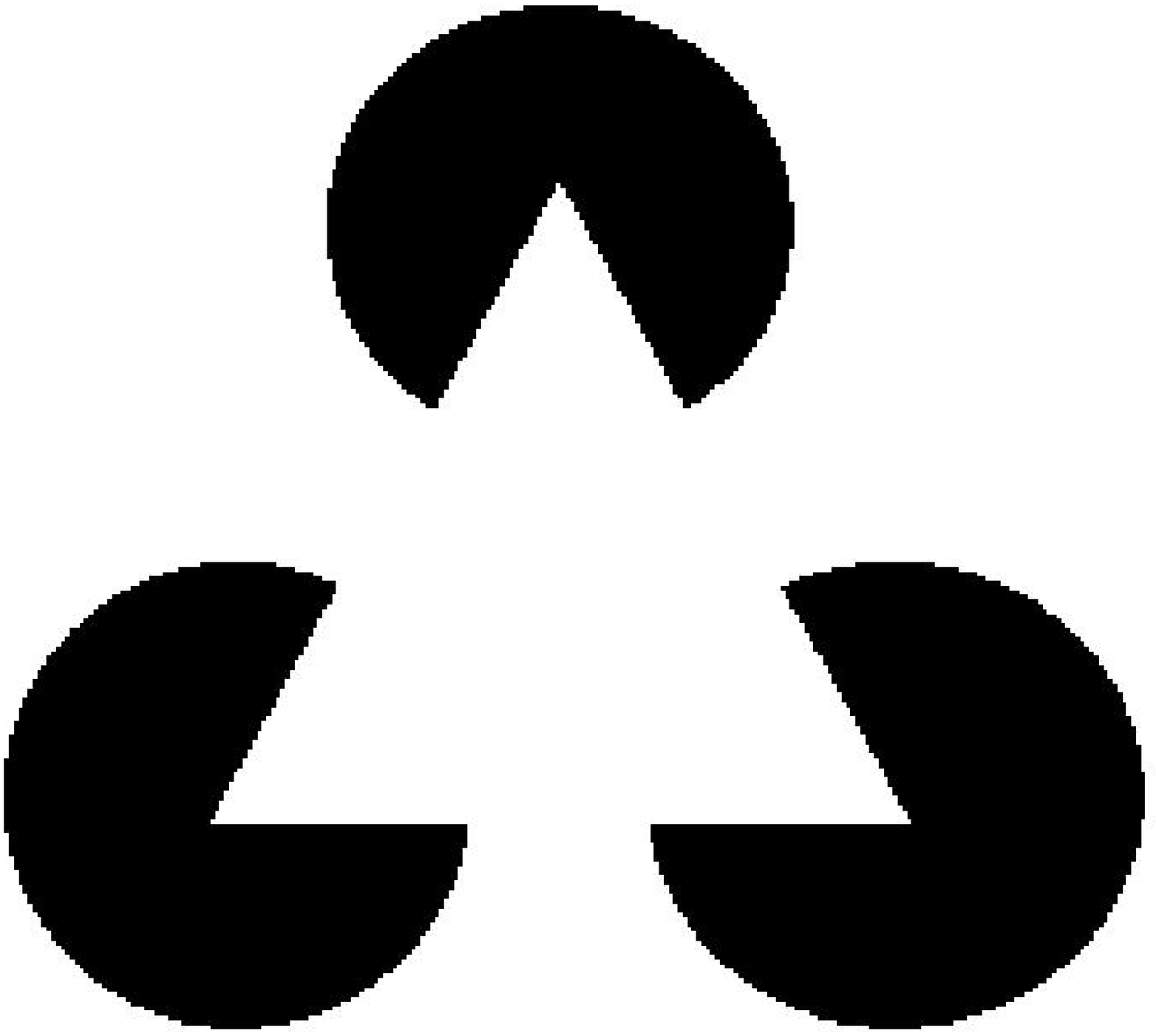}
  }
  \subfigure[Canyon Function $G$]{
   \includegraphics[width=3.5cm]{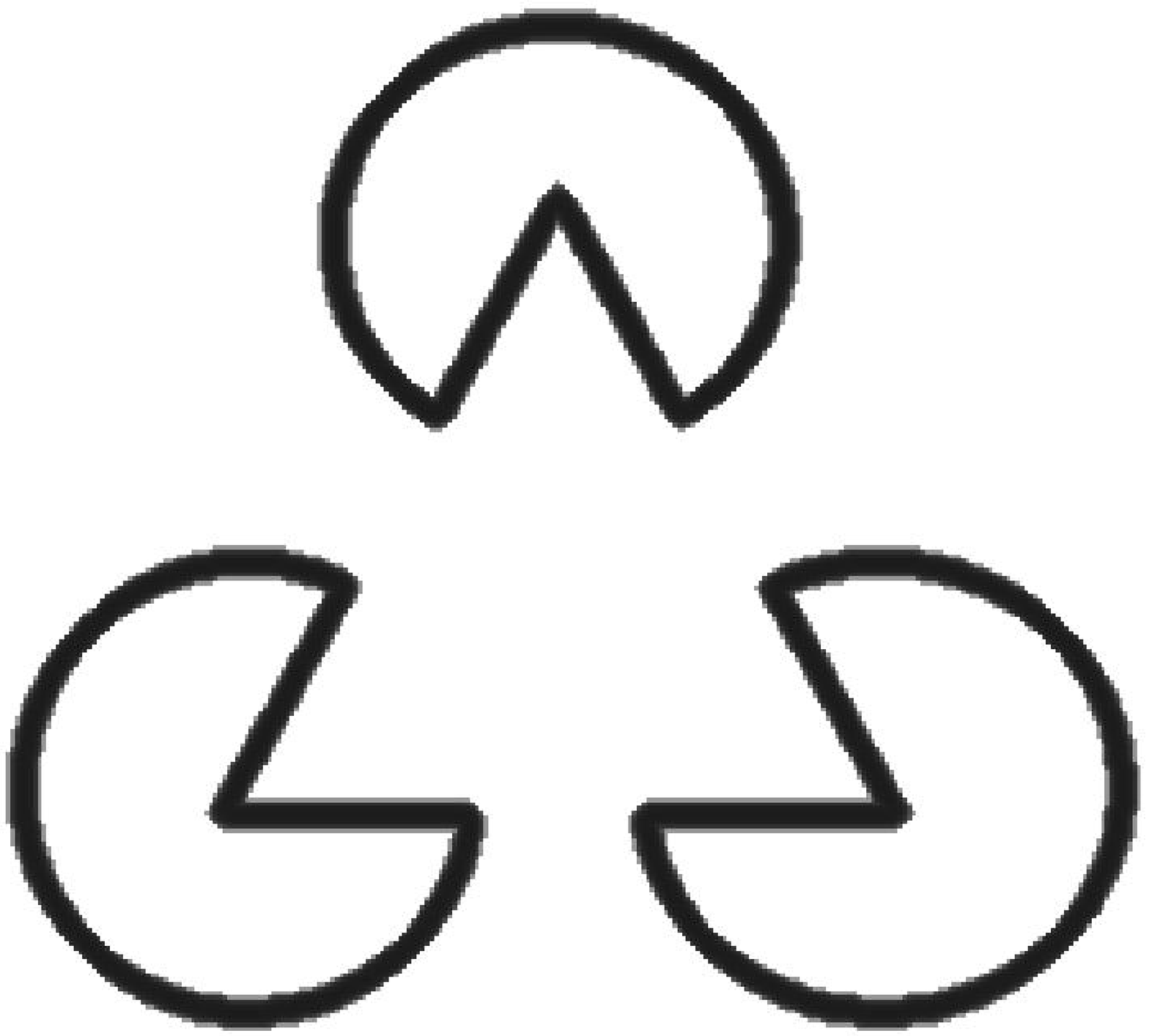}
  }
  \subfigure[Null Hypothesis $z_0(\bx)$]{
   \includegraphics[width=3.5cm]{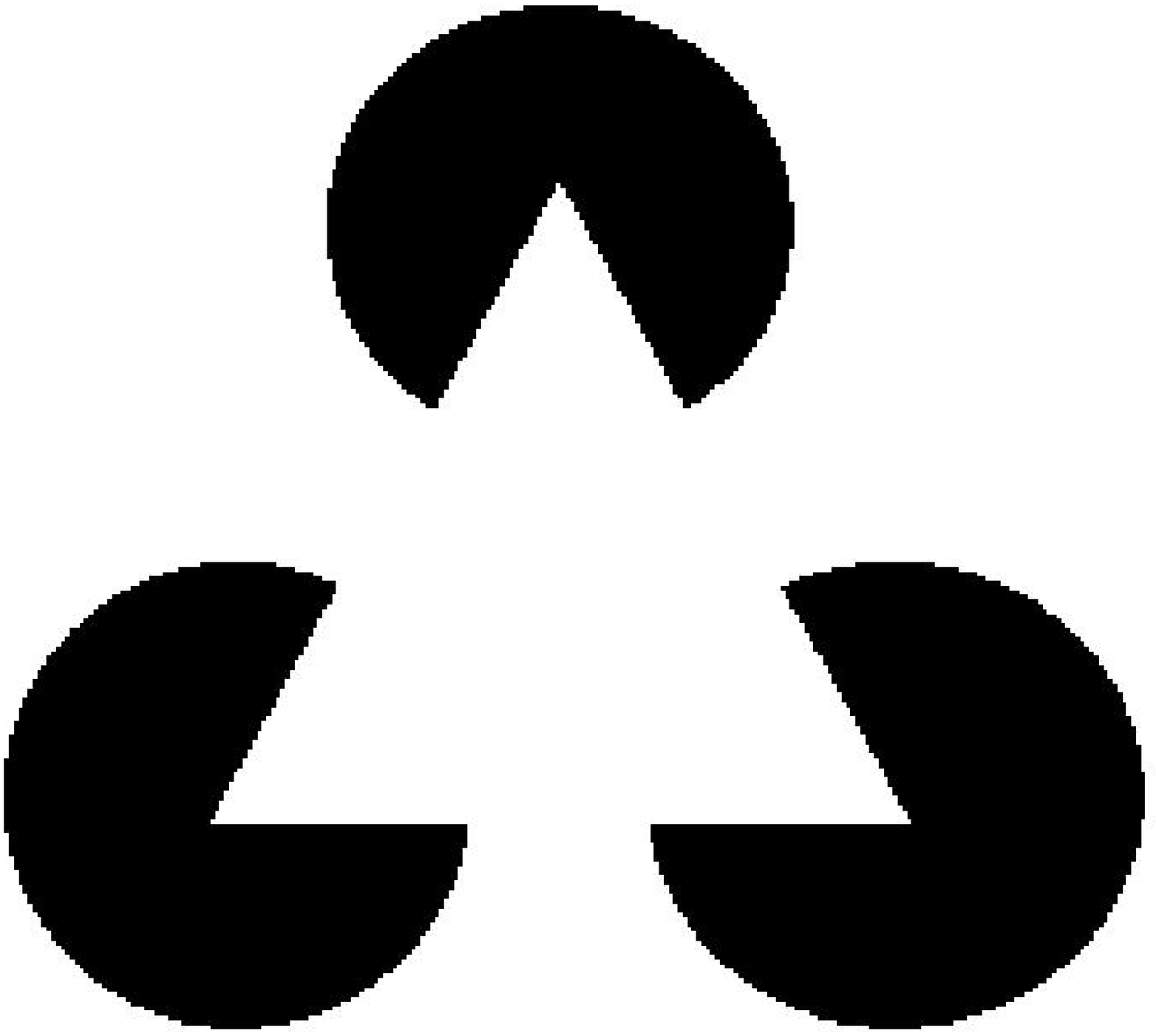}
  }\\
  \subfigure[$20$ Iterations]{
   \includegraphics[width=3.5cm]{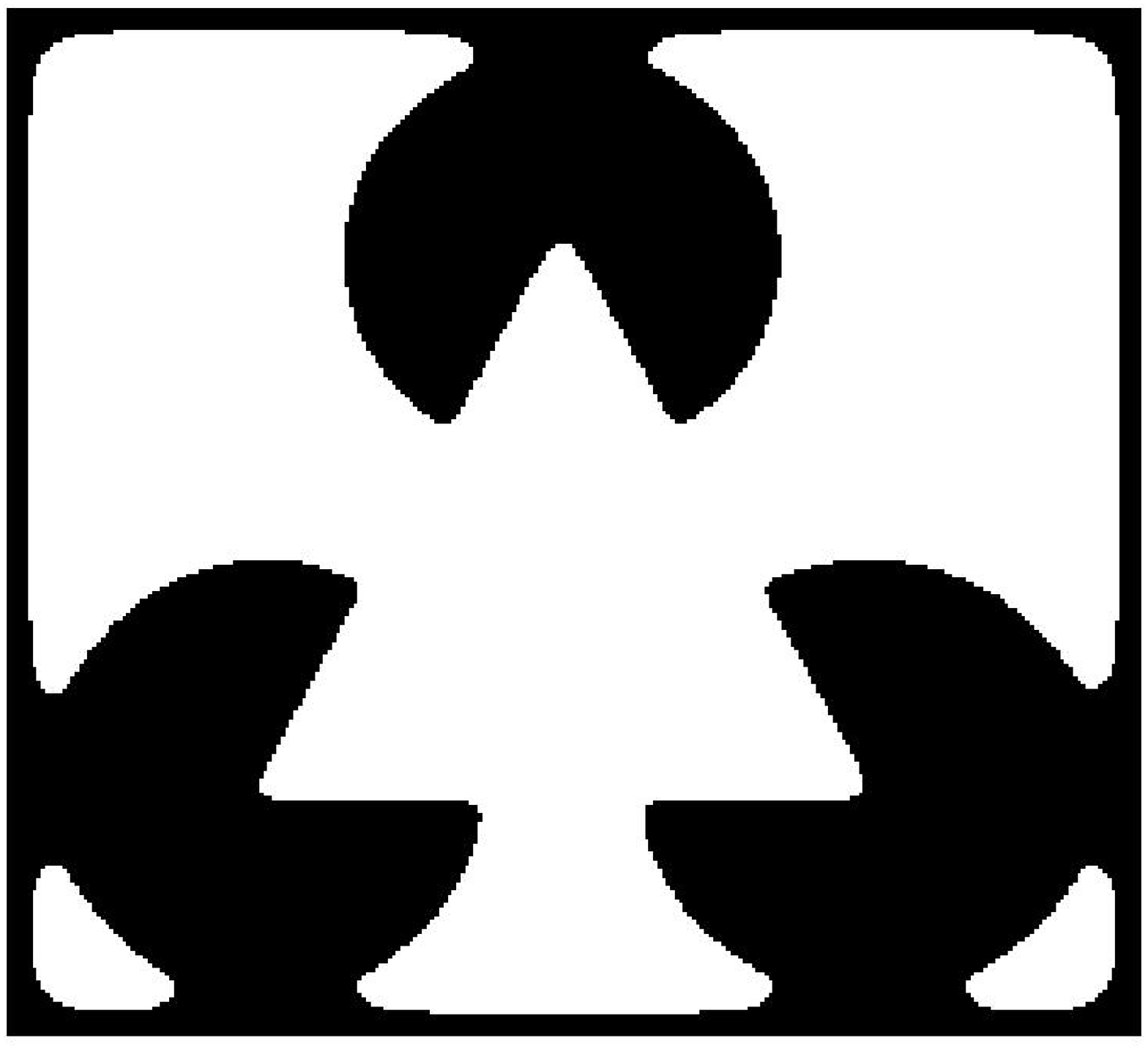}
  }
  \subfigure[$100$ Iterations]{
   \includegraphics[width=3.5cm]{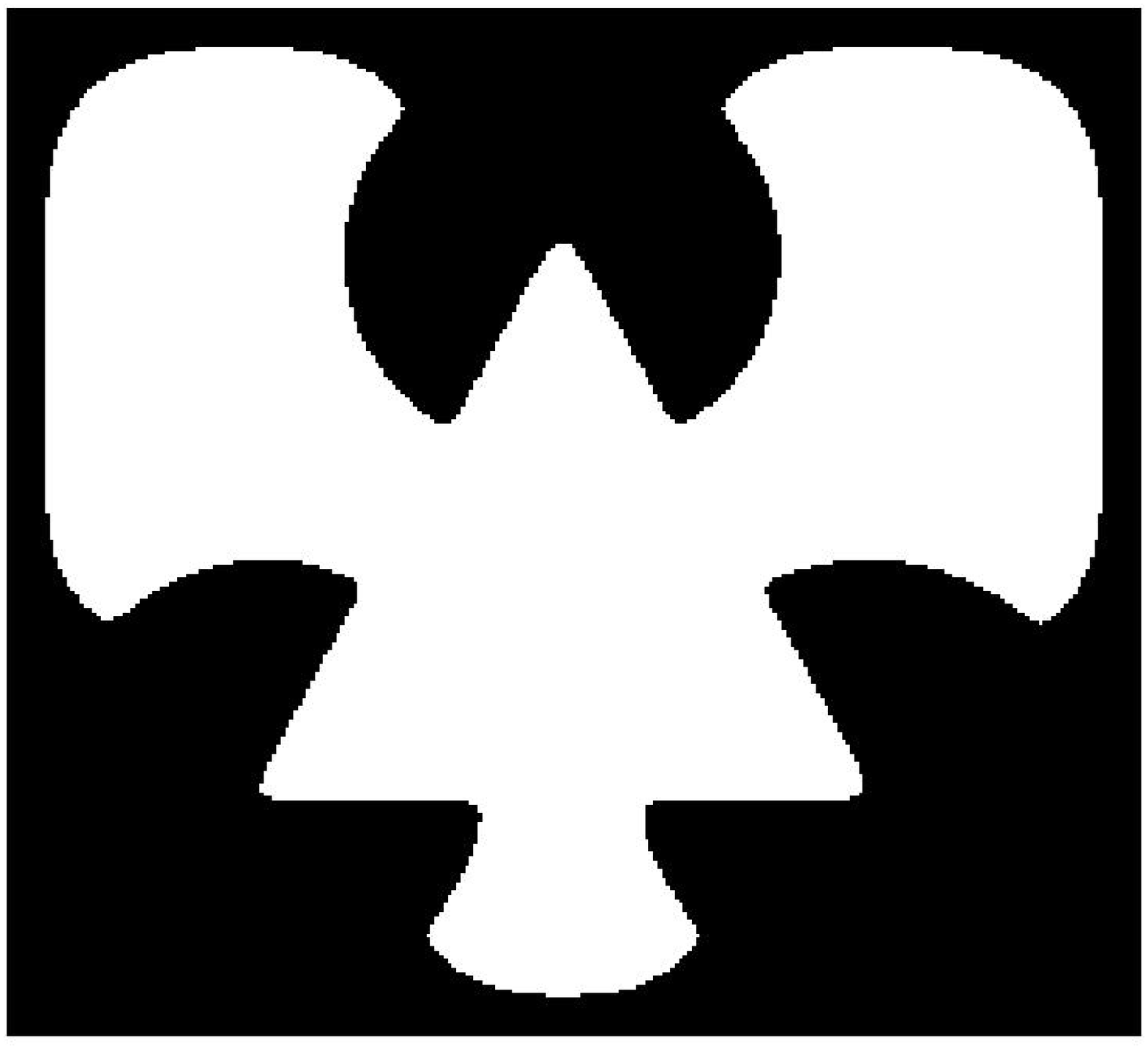}
  }
  \subfigure[Final Illusory Shape]{
   \includegraphics[width=3.5cm]{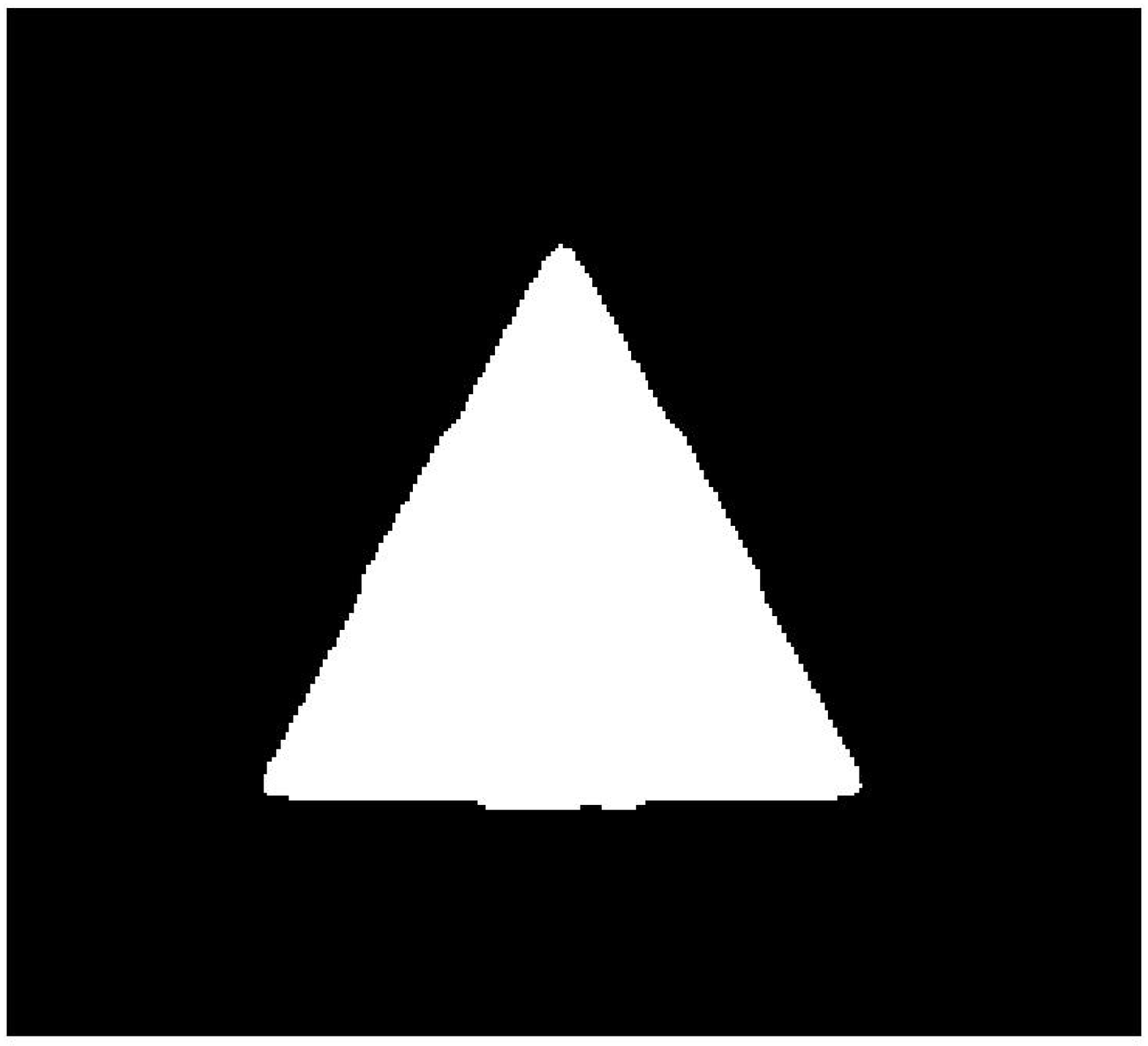}
  }
  \caption{The model and algorithm on the Kanizsa Triangle.}
  \label{fig:triangle}
\end{figure}
Fig.~\ref{fig:triangle} shows the numerical simulation on the classic example of Kanizsa's Triangle. In all the panels (b-f), the white colors represent phase values close to $1.0$ while the black ones to $0.0$. Panel (b) shows the canyon function that is fed into the model and algorithm. Panel (c) shows the null hypothesis $z_0(\bx)$ that conservatively assigns phase $1.0$ to all pixels outside $Q$. Panels (d) and (e) show two intermediate iterations before final convergence. In Panel (f) the final phase field is plotted after numerical convergence. It successfully captures the illusory triangle (up to the numerical precision).  Fig.~\ref{fig:circle} and Fig.~\ref{fig:butterfly_lady} demonstrate another two examples with more complex layouts.

\begin{figure}[h!]
  \centering
  \subfigure[Illusory Disk]{
   \includegraphics[width=3.7cm]{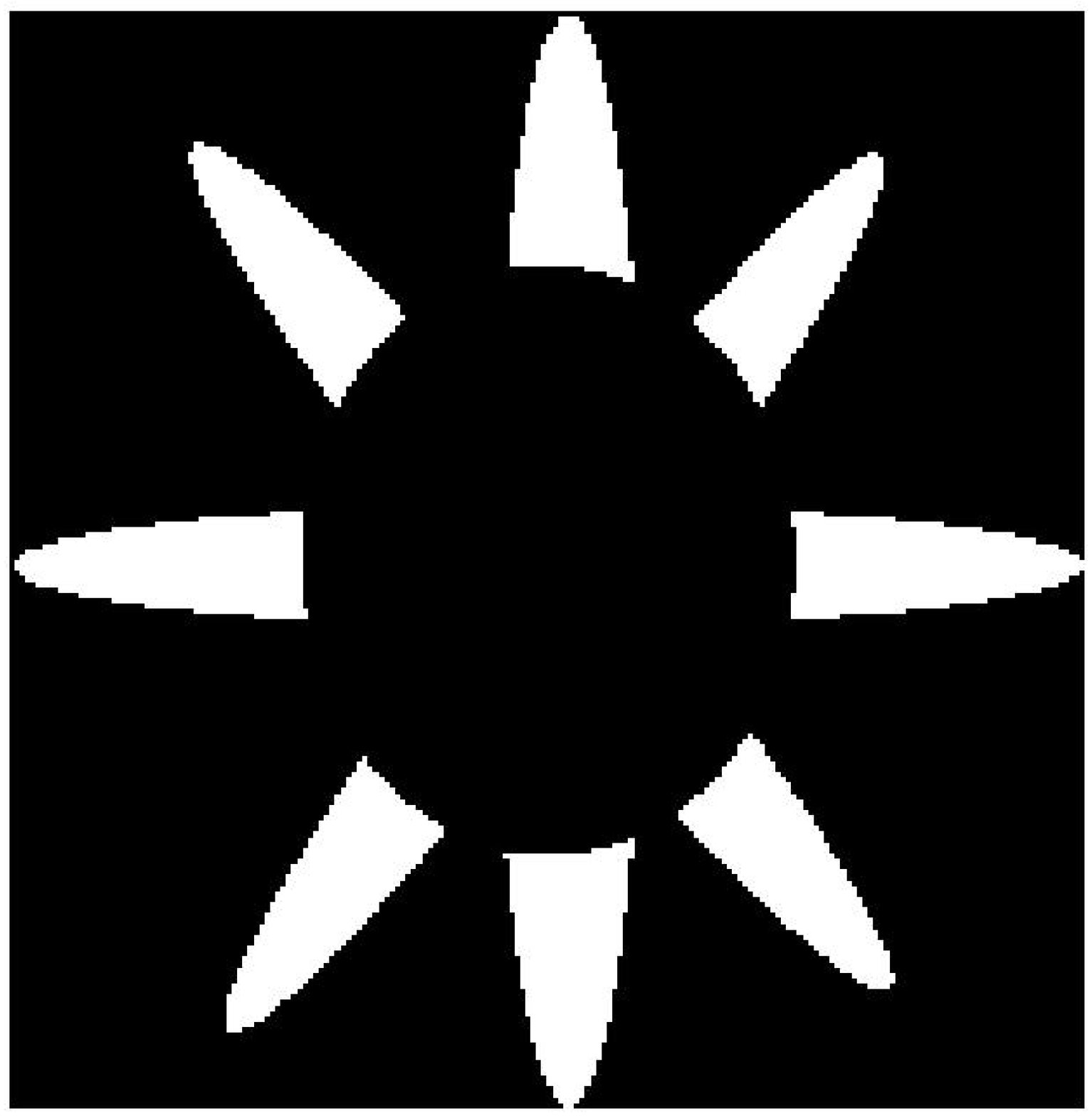}
  }
  \subfigure[Intermediate Snapshot]{
   \includegraphics[width=3.7cm]{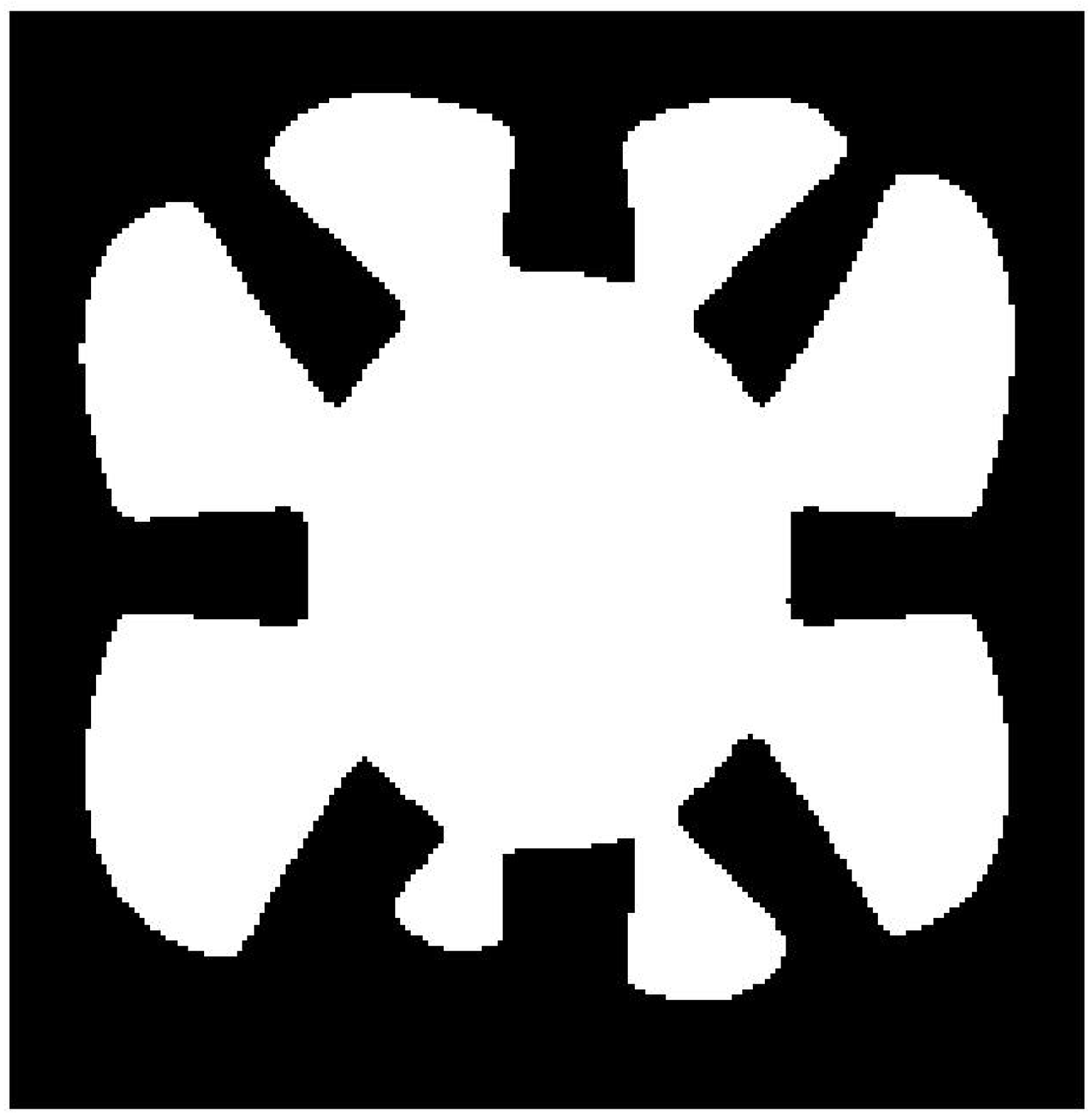}
  }
  \subfigure[Final Result]{
   \includegraphics[width=3.7cm]{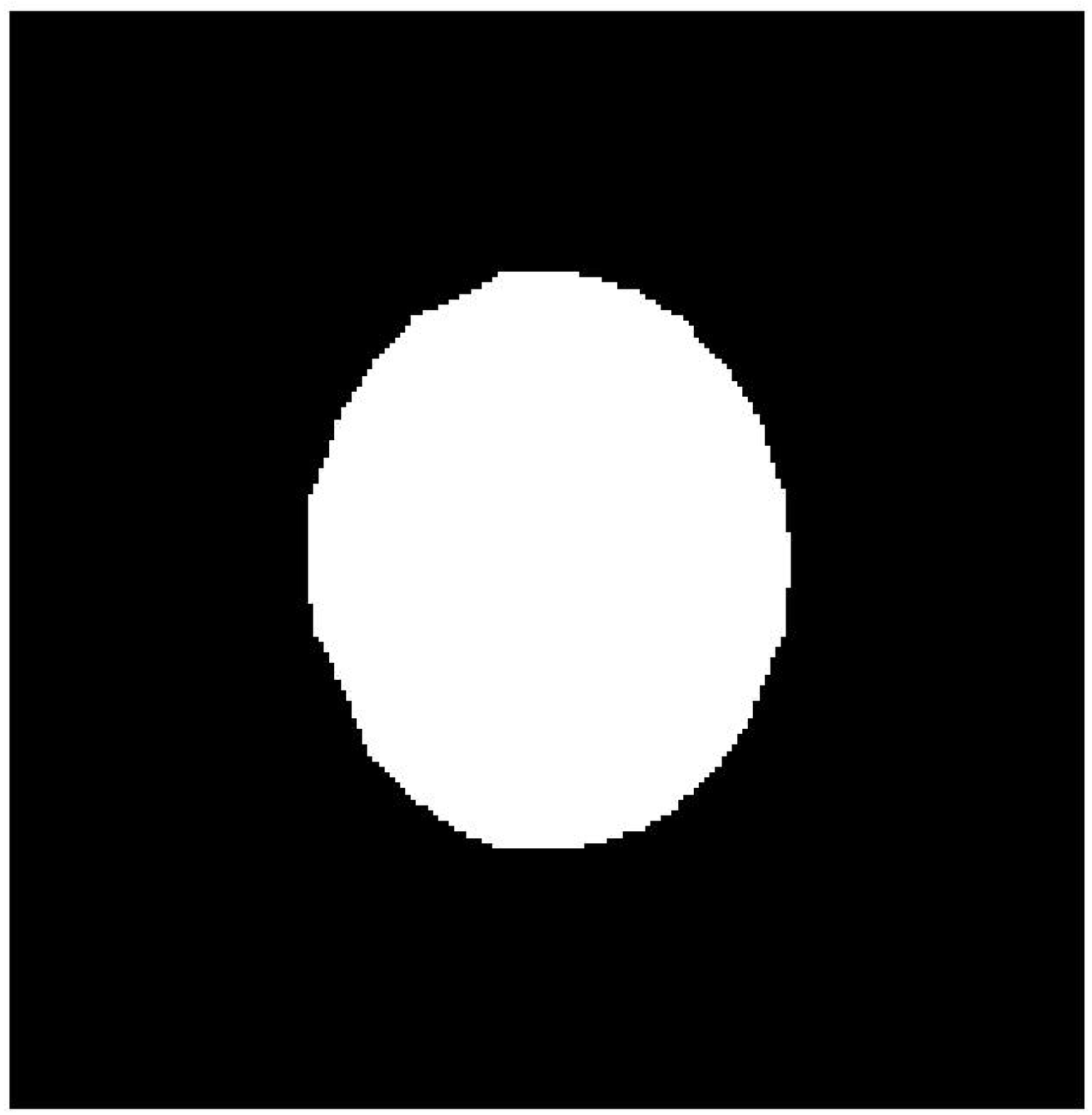}
  }
  \caption{The model and algorithm on an example with an illusory disk.}
  \label{fig:circle}
\end{figure}

\begin{figure}[h!]
  \centering
  \subfigure[Illusory Lady]{
   \includegraphics[height=4cm]{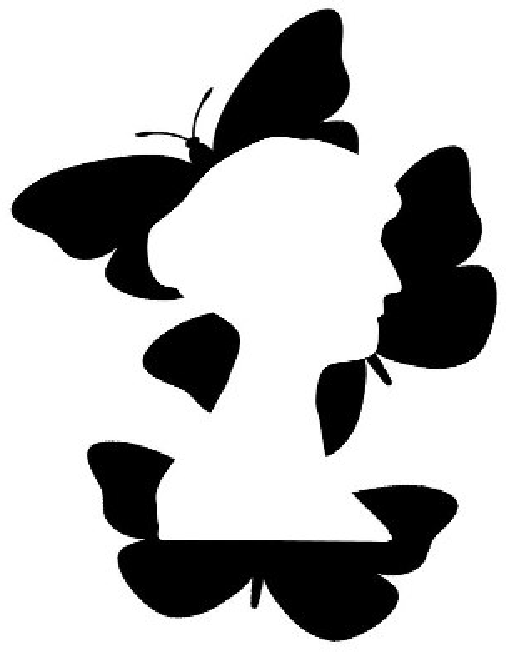}
  } \quad \quad \quad \quad
  \subfigure[Converged Result]{
   \includegraphics[height=4cm]{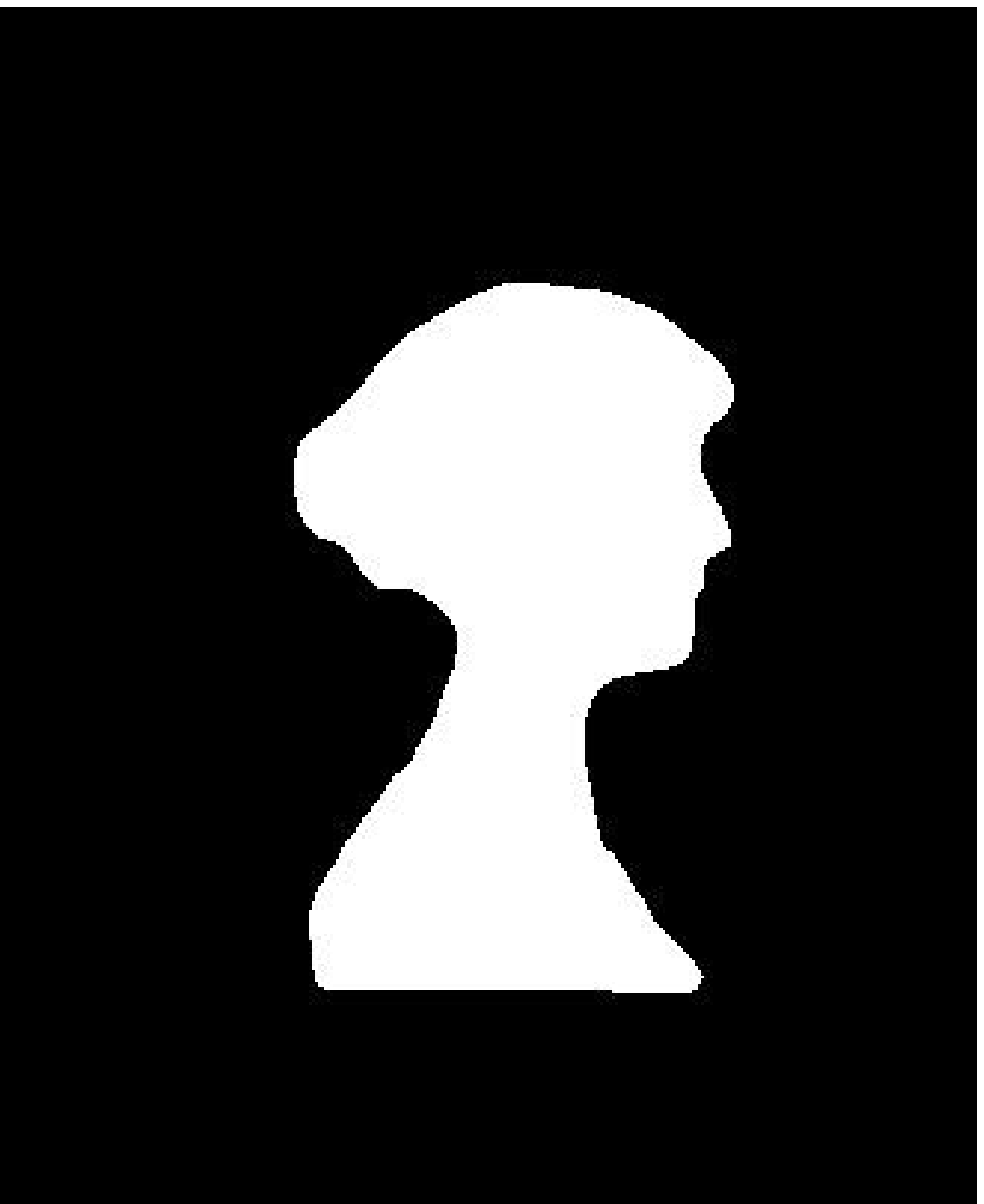}
  }
  \caption{The Butterfly Lady, with a more complex configuration $Q$. }
  \label{fig:butterfly_lady}
\end{figure}

Fig.~\ref{fig:twoobjs} shows an example consisting of two disjoint illusory shapes: an ellipse and a triangle. As a result, the phase field sequence $(z_n)$ is expected to experience a topological splitting operation during the iteration. Like the renowned level-set methodology of Osher and Sethian~\cite{oshset}, the phase-field approach is also very versatile in handling region merging or splitting.

\begin{figure}[h!]
  \centering
  \subfigure[Illusory Ellipse and Triangle]{
   \includegraphics[width=5cm]{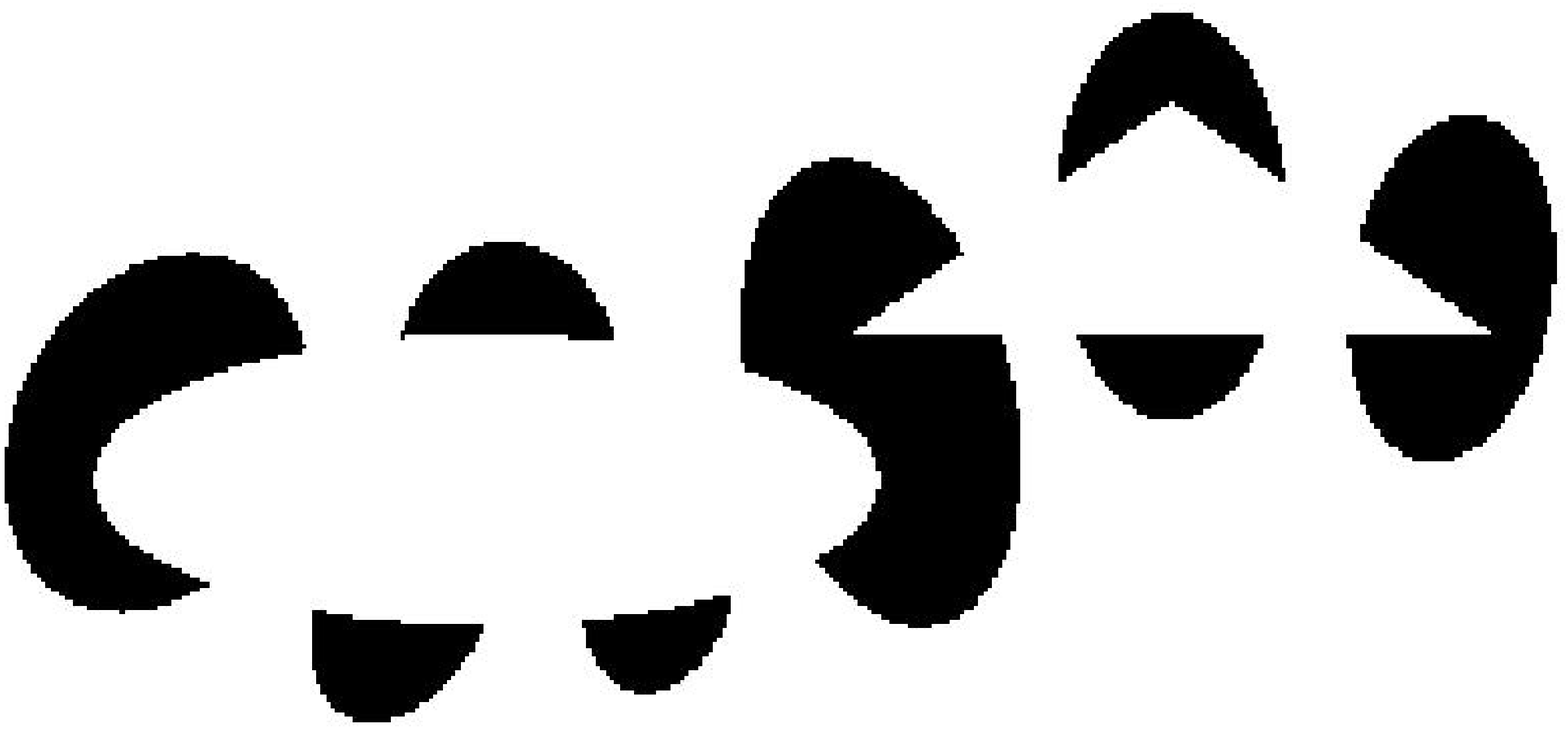}
  }
  \subfigure[Converged Result from the Model]{
   \includegraphics[width=5cm]{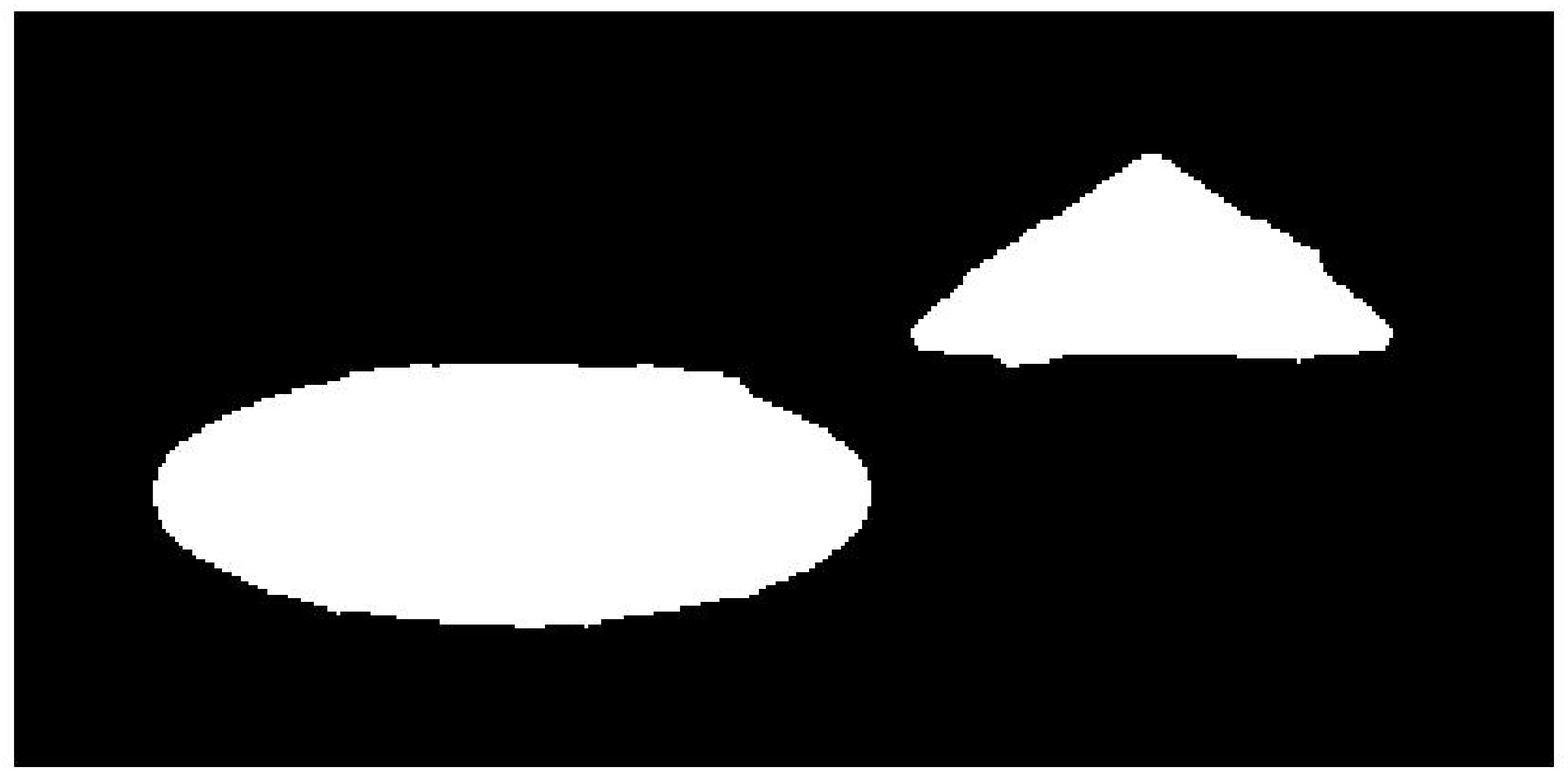}
  }
  \caption{Like the level-set methodology of Osher and Sethian~\cite{oshset}, the phase-field approach is versatile in handling topological changes like region splitting. }
  \label{fig:twoobjs}
\end{figure}

\section{Conclusion}
\label{sec:conclusion}
Inspired by the first-order variational illusory contour (VIC) model proposed in~\cite{illu_jungshen}, we have proposed a variational illusory shape (VIS) model based on the tool of phase transitions. The VIS model represents an illusory shape via phase values close to 1.0, and the rest by values close to 0.0. The phase transition is achieved by a variational energy formulated in the current work.

As for most non-quadratic phase transition models~\cite{gamma_March92,gamma_MarchDozio97,gamma_ShenPCMS05,gamma_ShenSoftMS06}, the proposed VIS model is non-convex. The zero field is the global optimum but uninteresting. To seek visually meaningful local optima, we have designed an iterative algorithm with a suitable initial guess, which could be considered as the {\em null hypothesis} in statistical testing. The null hypothesis assumes that there exists an illusory shape outside the given configuration. The iterative algorithm lets the pixels  ``vote'' collectively, until reaching the final consistent and stationary decision. Some key behaviors of the algorithm have been revealed through our analysis. And several generic numerical examples show the versatility of the proposed model and algorithm.

As in~\cite{illu_jungshen}, such lower-order models allow one to develop detailed analysis, but are necessarily limited in terms of applicability or performance. For example, illusory interpolation is often done via straight lines. Nevertheless, they help point towards more complex high-order models involving the curvature feature or Euler's elasticas~\cite{bImage_ChanShen,chakanshe,mum_elastica,illu2014_kangzhushen}, for example.

\begin{acknowledgements}
Jung has been supported by Basic Science Research Program through the National Research Foundation of Korea (NRF) of Korea (2012R1A1A1015492). Shen has been supported by the National Science Foundation (NSF) of USA.
\end{acknowledgements}


\end{document}